\documentclass[11pt]{article}

\usepackage{amsfonts}
\usepackage{latexsym}
\usepackage{color}
\usepackage{algorithmic}

\usepackage{hyperref}

\newtheorem{theorem}{Theorem}[section]
\newtheorem{lemma}[theorem]{Lemma}

\newtheorem{example}[theorem]{Example}

\newtheorem{statement}[theorem]{Statement}
\newtheorem{problem}[theorem]{Problem}

\newcommand{\bff}{\mbox{$\mbox{\boldmath $f$}$}}%added by ph
\newcommand{\bfe}{\mbox{$\mbox{\boldmath $e$}$}}%added by ph
\newcommand{\bfg}{\mbox{$\mbox{\boldmath $g$}$}} %added by ph
\newcommand{\bfh}{\mbox{$\mbox{\boldmath $h$}$}} %added by ph
\newcommand{\bfk}{\mbox{$\mbox{\boldmath $k$}$}} %added by ph
\newcommand{\bfr}{\mbox{$\mbox{\boldmath $r$}$}} %added by ph
\newcommand{\bfv}{\mbox{$\mbox{\boldmath $v$}$}} %added by ph
\newcommand{\bfw}{\mbox{$\mbox{\boldmath $w$}$}} %added by ph
\newcommand{\bfx}{\mbox{$\mbox{\boldmath $x$}$}} %added by ph
\newcommand{\bfy}{\mbox{$\mbox{\boldmath $y$}$}} %added by ph
\newcommand{\bfz}{\mbox{$\mbox{\boldmath $z$}$}} %added by ph

\newcommand{\sbfk}{\mbox{\boldmath \scriptsize $k$}} %added by ph
\newcommand{\sbfx}{\mbox{\boldmath \scriptsize $x$}}%added by ph

\newcommand{\bxi}{\mbox{$\mbox{\boldmath $\xi$}$}} %added by ph
\newcommand{\bfeta}{\mbox{$\mbox{\boldmath $\eta$}$}} %added by ph
\newcommand{\bfzero}{\mbox{\boldmath $0$}}

\begin{document}

\title{Optimal splitting of Parseval frames using Walsh matrices}

\author{Amie Albrecht, Phil Howlett, Geetika Verma \\
\\
Scheduling and Control Group, \\
Centre for Industrial and Applied Mathematics, \\
School of Information Technology and Mathematical Sciences, \\
University of South Australia, \\
Mawson Lakes, Australia, 5095. \\
Email: phil.howlett@unisa.edu.au.}

\maketitle

\maketitle

\begin{abstract}
\noindent In 2014 Adam Marcus, Daniel Spielman and Nikhil Srivastava used random vectors to prove a key discrepancy theorem and in so doing gave a positive answer to the long-standing Kadison\textendash Singer Problem.  In this paper we use Walsh matrices to construct a class of natural frames in Euclidean space and discuss how these frames relate to the key discrepancy theorem.
\end{abstract}



\section{Introduction}
\label{sec:intro}

In 1959 Richard Kadison and Isadore Singer \cite{kad1} formulated a problem in quantum mechanics that later became one of the iconic mathematical questions of the twentieth century.  The problem is now known as the Kadison\textendash Singer Problem (KSP).

\begin{problem}[KSP]
\label{p1}
Does every pure state on the algebra of bounded diagonal operators acting on the Hilbert space of square summable complex-valued sequences have a unique extension to a regular state on the algebra of all bounded operators? $\hfill \Box$
\end{problem}

Following a finite-dimensional reformulation \cite{and1} by Joel Anderson in 1979 and further reduction to an equivalent problem in discrepancy theory \cite{wea1} by Nik Weaver in 2004, a positive answer to KSP was eventually found \cite{mar1,mar2} by Adam Marcus, Daniel Spielman and Nikhil Srivastava in 2013.  We will not attempt a detailed explanation of KSP but instead refer readers to the excellent review article \cite{har1} by Nick Harvey. The  Marcus\textendash Spielman\textendash Srivistava Discrepancy Theorem (MSSDT) was a basic platform for the eventual solution of KSP and is a central theme in our paper.

\begin{theorem}[MSSDT]
\label{ksdiscrepancy}
If $\bfv_1,\ldots,\bfv_n \in {\mathbb C}^m$ are such that $\|\bfv_j\|^2 \leq \alpha$ for all $j=1,\ldots,n$ and
\begin{equation}
\label{ntfmatrixsum}
\sum_{j=1}^n \bfv_j \bfv_j^* = I
\end{equation}
where $I \in {\mathbb C}^{m \times m}$ is the unit matrix then there is a partition of the index set ${\mathcal J} = \{1,\ldots,n\} \subset {\mathbb N}$ into two disjoint subsets ${\mathcal J}_1$ and ${\mathcal J}_2$ such that
\begin{equation}
\label{ksdcon}
\left \| \sum_{j\, \in\, {\mathcal J}_k} \bfv_j {\bfv_j}^* \right \|_2 \leq \left( \frac{1}{\sqrt{2}} + \sqrt{\alpha} \right)^2
\end{equation}
for each $k=1,2$.  The norm used here is the $2$-norm.  Note that $(\ref{ntfmatrixsum})$ implies $m \leq n$.  $\hfill \Box$
\end{theorem}

MSSDT was a fundamental stepping stone in the ultimately successful quest \cite{mar1, mar2} for a positive answer to  KSP \cite{kad1}.  In fact MSSDT also implies the truth of the Weaver Discrepancy Statement (WDS) proposed earlier in 2004 by Nik Weaver \cite{wea1} as a mechanism for solving KSP.

\begin{statement}[WDS]
\label{wds}
Let $\bfv_1,\ldots,\bfv_n \in {\mathbb C}^m$ satisfy $\|\bfv_j\|^2 \leq \alpha$ for each $j=1,\ldots,n$ and suppose that
\begin{equation}
\label{ntfquadform}
\sum_{j=1}^n |{\bfv_j}^*\bfx|^2 = 1
\end{equation}
for all $\bfx \in {\mathbb C}^m$ with $\|\bfx\|=1$.  Then we can partition ${\mathcal J} = \{1,\ldots,n\}$ into two disjoint sets ${\mathcal J}_1, {\mathcal J}_2$ such that
\begin{equation}
\label{wdscon}
\left |  \sum_{j \in {\mathcal J}_k} |{\bfv_j}^* \bfx|^2 - \frac{1}{2} \right | \leq 5 \sqrt{\alpha}
\end{equation}
for each $k=1,2$ and all $\bfx \in {\mathbb C}^m$ with $\|\bfx\|=1$.  If we write ${\mathcal J}_1 = \{h(1),\dots,h(p)\}$ where $1 \leq h(1) < \cdots < h(p) \leq n$ and ${\mathcal J}_2 = \{k(1),\ldots,k(q)\}$ where $1 \leq k(1) < \cdots < k(q) \leq n$ with $p+q=n$ and define $V_1 = [\bfv_{h(1)},\ldots,\bfv_{h(p)}]$ and $V_2 = [\bfv_{k(1)},\ldots,\bfv_{k(p)}]$ then $(\ref{wdscon})$ can be rewritten as
\begin{equation} 
\label{wdsmatrix}
\left \| \rule{0cm}{0.4cm} \hspace{1mm} V_kV_k^* - \frac{I}{2} \hspace{1mm} \right \|_2 \leq 5 \sqrt{\alpha}, 
\end{equation}
where $I \in {\mathbb C}^{m \times m}$ is the unit matrix, for each $k=1,2$.  $\hfill \Box$
\end{statement}

An advantage of  WDS is that it allows us to interpret the discrepancy in terms of quadratic forms.  WDS says that any quadratic form expressed as a sum of small rank one quadratic forms can be split into two almost equal parts.  The {\em discrepancy} is represented as the difference in values on the surface of the unit sphere of the two constituent quadratic forms with each one expressed as a sum of small rank one quadratic forms.  A key motivation for our paper is the close connection between MSSDT and the theory of frames \cite{cas2,chr1} in finite dimensional Euclidean space.
 
\subsection{Motivation}
\label{m} 
 
Conditions (\ref{ntfmatrixsum}) and (\ref{ntfquadform}) are equivalent.  A set of vectors $\{\bfv_1,\ldots,\bfv_n\} \in {\mathbb C}^m$ that satisfies these conditions is said to form a {\em Parseval frame} or {\em normalized tight frame} in ${\mathbb C}^m$.  Based on MSSDT one of our key motivations was to explore the connection between discrepancy theory and Parseval frames in finite-dimensional Euclidean space.  Our second motivation is less obvious.  In a recent note \cite{sri1}, Nikhil Srivastava wrote that in general $\ldots$ the presence of large vectors (in a frame) {\em is an obstruction to the existence of a low discrepancy partition}.  Thus we decided to investigate Parseval frames in which all vectors are the same size.  It is known that Parseval frames in finite-dimensional Euclidean spaces are closely related to orthogonal matrices.  The Walsh matrices are a well-known collection of real symmetric matrices where all elements have magnitude $1$ and the sets of row and column vectors are each sets of mutually orthogonal vectors.  Thus we were led to a discussion of discrepancy theory for Parseval frames defined by Walsh matrices. 

\subsection{Tight frames in finite dimensional Euclidean space}
\label{ntffdes}

If the set of vectors $\{\bfv_1,\ldots,\bfv_n\} \in {\mathbb C}^m$ satisfies the condition
\begin{equation}
\label{tf}
\sum_{j=1}^n |\bfv_j^*\bfx|^2 = c
\end{equation}
for some $c > 0$ and all $\bfx \in {\mathbb C}^m$ with $\|\bfx\| = 1$ then $\{\bfv_1,\ldots,\bfv_n\}$ is said to form a {\em tight frame} in ${\mathbb C}^m$ with frame constant $c$.  In such cases we must have $m \leq n$.  If we define the {\em pre-frame operator} $V = [\bfv_1,\ldots,\bfv_n] \in {\mathbb C}^{m \times n}$ then
$$
\sum_{j=1}^n |{\bfv_j}^* \bfx|^2 = \bfx^* \left[ \sum_{j=1}^n \bfv_j {\bfv_j}^* \right] \bfx = \bfx^* VV^* \bfx = \bfx^* S \bfx
$$
where $S = VV^* \in {\mathbb C}^{m \times m}$ is called the {\em frame operator}.  The frame operator is self-adjoint and positive.  Hence it is invertible.  For all $\bfx \in {\mathbb C}^m$ we can write
$$
\bfx = \sum_{j=1}^n \left( {\bfv_j}^* S^{-1} \bfx \right) \bfv_j 
$$
where the coefficients ${\bfv_j}^* S^{-1} \bfx$ for each $j=1,\ldots,m$ are called the frame coefficients for $\bfx$.  In the case where $c=1$ the condition (\ref{tf}) reduces to (\ref{ntfquadform}) and the frame becomes a Parseval frame.  Condition (\ref{ntfquadform}) can now be rewritten as  
\begin{equation}
\label{ntfframeop}
\bfx^*(S - I_m)\bfx = 0
\end{equation}
for all $\bfx \in {\mathbb C}^m$ with $\bfx \neq \bfzero$.  Thus for a Parseval frame we must have $S = I_m$ and the frame representation reduces to
\begin{equation}
\label{ntfrep}
\bfx = \sum_{j=1}^n \left( {\bfv_j}^* \bfx \right) \bfv_j 
\end{equation}
for all $\bfx \in {\mathbb C}^m$.  This formula suggests another more familiar formula.  If we define $W = V^*$ and write $W = [\bfw_1,\ldots,\bfw_m] \in {\mathbb C}^{n \times m}$ then the condition $S = VV^* = I_m$ can be rewritten as $W^*W = I_m$ and this means that the set $\{\bfw_1,\ldots,\bfw_m\}$ forms an orthonormal set in ${\mathbb C}^n$.  If we extend this set to an orthonormal basis $\{\bfw_1,\ldots,\bfw_n\}$ and write $H = [H_1 \mid H_2] = [\bfw_1, \ldots, \bfw_m \mid \bfw_{m+1},\ldots,\bfw_n] \in {\mathbb C}^{n \times n}$ then $H_1 = W$.  Let $G = H^*$ and write
$$
G = \left[ \begin{array}{cc}
G_1 \\ \hline
G_2 \end{array} \right]
$$
where $G_1 = {H_1}^* = V$.  Thus we may write
$$
G = \left[ \begin{array}{cccc}
\bfv_1 & \bfv_2 & \cdots & \bfv_n \\ \hline
\bfr_1 & \bfr_2 & \cdots & \bfr_n \end{array} \right]
$$
where $\{ \bfr_1,\ldots,\bfr_n\} \in {\mathbb C}^{(n-m) \times n}$.  The matrix $G = [\bfg_1,\ldots,\bfg_n]  \in {\mathbb C}^{n \times n}$ is orthogonal and the set of vectors $\{\bfg_1,\ldots,\bfg_n\} \in {\mathbb C}^n$ forms an orthonormal basis for ${\mathbb C}^n$.  Since
$$
\bfg_j = \left[ \begin{array}{cc}
\bfv_j \\
\bfr_j \end{array} \right]
$$
for each $j=1,\ldots,n$ the standard representation for a vector
$$
\bfz = \left[ \begin{array}{cc}
\bfx \\
\bfr \end{array} \right] \in {\mathbb C}^n
$$
in terms of the orthonormal basis $\{\bfg_1,\ldots \bfg_n\}$ is given by
\begin{equation}
\label{obrep}
\bfz = \sum_{j=1}^n \left( {\bfg_j}^* \bfz \right) \bfg_j \iff \left[ \begin{array}{cc}
\bfx \\
\bfr \end{array} \right] = \sum_{j=1}^n \left( {\bfv_j}^* \bfx + {\bfr_j}^*\bfr \right) \left[ \begin{array}{cc}
\bfv_j \\
\bfr_j \end{array} \right].
\end{equation}
For vectors in the subspace $S_m \subseteq {\mathbb C}^n$ defined by $\bfr = \bfzero$ the representation in (\ref{obrep}) reduces to
\begin{equation}
\label{efrep}
\left[ \begin{array}{cc}
\bfx \\
\bfzero \end{array} \right] = \sum_{j=1}^n \left( {\bfv_j}^* \bfx \right) \left[ \begin{array}{cc}
\bfv_j \\
\bfr_j \end{array} \right]
\end{equation}
which is essentially the same representation as (\ref{ntfrep}).  It follows that $\sum_{j=1}^n ({\bfv_j}^* \bfx) \bfr_j = \bfzero$.  In fact we can see that if $\bfz \in S_m$ then
$$
\sum_{j=1}^n | {\bfg_j}^* \bfz |^2 = \sum_{j=1}^n \left( | {\bfv_j}^* \bfx |^2 + | {\bfr_j}^* \bfzero|^2 \right) = \sum_{j=1}^n | {\bfv_j}^* \bfx |^2 = 1
$$
for all $\bfz \in S_m$ with $\|\bfz\|^2 = 1$.  Hence the set $\{ \bfg_1,\ldots,\bfg_n\} \in {\mathbb C}^n$ defines a Parseval frame for the $m$-dimensional subspace $S_m \subseteq {\mathbb C}^n$.  The frame for $S_m$ defined by $G$ is simply the original frame defined by $V$ embedded into ${\mathbb C}^n$.  We could write $G \cong V$.

For a Parseval frame defined by a pre-frame operator $V$ the frame operator $S = VV^* = W^*W = I_m \in {\mathcal B}({\mathbb C}^m)$ is the identity mapping.  The related operator $R = V^*V = WW^* = J_n \in {\mathcal B}({\mathbb C}^n)$ satisfies the equations $J_n^2 = J_n$ and  $J_n W = W$.  Thus $J_n$ is a projection onto the column space of $W$.

\subsection{The Walsh functions and Walsh matrices}
\label{wm}

The {\em Walsh functions} $W_k:[0,1] \rightarrow \{-1,1\}$ for $k \in {\mathbb N} - 1$ can be defined as follows.  Choose $m \in {\mathbb N}$ and let each $k < 2^m$ be represented in binary form as
$$
k = k_m \cdots k_1 \iff k = \sum_{s=1}^m k_s\, 2^{s-1} \iff \bfk = (k_1,\ldots,k_m,0,0,\ldots) \in \{0,1\}^{\infty}
$$
and let each $x \in [0,1]$ be represented in binary form as
$$
x = 0.x_1x_2 \cdots \iff x = \sum_{s=1}^{\infty} x_s 2^{-s} \iff \bfx = (x_1,x_2,\ldots) \in \{0,1\}^{\infty}
$$
where no expansion is permitted with $x_s = 1$ for all $s \geq n$ for some $n = n(x) \in {\mathbb N}$.  Then we have
$$
W_k(x) = (-1)^{p(\sbfk,\sbfx)}
$$
for each $k < 2^m$ and each $x \in [0,1]$ where $p(\bfk,\bfx) = \sum_{s=1}^m k_s x_s$.  The Walsh functions form a complete orthonormal set in the Hilbert space $L^2[0,1]$.  They were introduced in a 1923 paper \cite{wal1} by Joseph Walsh and have since found wide application in digital signal processing.  In this regard a fundamental requirement was the development of efficient computation routines for Walsh matrices and the associated function representations using Walsh series and Walsh transforms.  For a detailed account see \cite{sch1}.  See also \cite{fin1} and references therein and some recent work on the construction of wavelet frames in Walsh analysis \cite{far1, far2, far3}.  Importantly we note that Per Enflo used Walsh series to prove a celebrated result \cite{enf1} that there exist separable Banach spaces with no Schauder basis.

The Walsh functions are closely related to the Walsh matrices which are our particular interest in this paper.  Let $n = 2^r$ for some $r \in {\mathbb N}$.  The {\em Walsh matrix} $Y = Y_r \in {\mathbb C}^{n \times n}$ can be efficiently computed using the recursive Sylvester construction defined by the M{\sc atlab} algorithm \medskip

\begin{algorithmic}
\STATE $Y = [1]$;
\STATE {\bf for} {$i = 1: r$}
\STATE $Y = [Y,\hspace{0.2cm}  Y; Y,\hspace{0.2cm} -Y]$;
\STATE {\bf end}
\end{algorithmic} 

\noindent The matrix $Y$ is real symmetric with $y_{ij} = \pm 1$ for all $i, j \in \{1,\ldots,n\}$ and  $Y^*Y = nI$.  The columns $\{\bfy_1,\ldots,\bfy_n\}$ form an orthogonal basis for ${\mathbb C}^n$ with $\|\bfy_j\| = \sqrt{n}$ for all $j=1,\ldots,n$.  If we choose $m \leq n$ and define $W = [\bfy_1,\ldots,\bfy_m]/\sqrt{n} \in {\mathbb C}^{n \times m}$ and $V = W^*$ and if we write $V = [\bfv_1,\ldots,\bfv_n] \in {\mathbb C}^{m \times n}$ then the columns of $V$ define a Parseval frame in ${\mathbb C}^m$ which consists of $n = 2^r$ vectors $\bfv_1,\ldots,\bfv_n$ with $\| \bfv_1 \|^2 = \cdots = \| \bfv_n \|^2 = m/n$.

We have
$$
Y_1 = \left[ \begin{array}{r|r}
1 & 1 \\ \hline
1 & -1\end{array} \right],
$$
$$
Y_2 = \left[ \begin{array}{rr|rr}
1 & 1 & 1 & 1 \\
1 & -1 & 1 & -1 \\ \hline
1 & 1 & -1 & -1 \\
1 & -1 & -1 & 1 \end{array} \right]
$$
and
$$
Y_3 = \left[ \begin{array}{rrrr|rrrr}
1 & 1 & 1 & 1 & 1 & 1 & 1 & 1 \\
1 & -1 & 1 & -1 & 1 & -1 & 1 & -1 \\
1 & 1 & -1 & -1 & 1 & 1 & -1 & -1 \\
1 & -1 & -1 & 1 & 1 & -1 & -1 & 1 \\ \hline
1 & 1 & 1 & 1 & -1 & -1 & -1 & -1 \\
1 & -1 & 1 & -1 & -1 & 1 & -1 & 1 \\
1 & 1 & -1 & -1 & -1 & -1 & 1 & 1 \\
1 & -1 & -1 & 1 & -1 & 1 & 1 & -1 \end{array} \right].
$$
Note that these matrices are known as the {\em Walsh} matrices \cite{fin1} using the {\em natural ordering} and they are a special case of the {\em Hadamard} matrices \cite{hed1}.  The construction described in the M{\sc atlab} algorithm is due to Sylvester \cite[Section 3.1]{hed1}.  The Walsh matrices can be presented with various different orderings of the rows and columns.  The {\em sequency ordering} \cite{fin1} is defined by ordering the rows according to the number of sign changes in each row.  Thus, with this ordering, we have
$$
Z_1 = \left[ \begin{array}{rr}
1 & 1 \\ 
1 & -1\end{array} \right],
$$
$$
Z_2 = \left[ \begin{array}{rrrr}
1 & 1 & 1 & 1 \\
1 & 1 & -1 & -1 \\
1 & -1 & -1 & 1 \\ 
1 & -1 & 1 & -1\end{array} \right]
$$
and
$$
Z_3 = \left[ \begin{array}{rrrrrrrr}
1 & 1 & 1 & 1 & 1 & 1 & 1 & 1 \\
1 & 1 & 1 & 1 & -1 & -1 & -1 & -1 \\
1 & 1 & -1 & -1 & -1 & -1 & 1 & 1 \\
1 & 1 & -1 & -1 & 1 & 1 & -1 & -1 \\
1 & -1 & -1 & 1 & 1 & -1 & -1 & 1 \\
1 & -1 & -1 & 1 & -1 & 1 & 1 & -1 \\
1 & -1 & 1 & -1 & -1 & 1 & -1 & 1 \\
1 & -1 & 1 & -1 & 1 & -1 & 1 & -1 \end{array} \right].
$$
The advantage of the sequency ordering is that row $k+1$ of $Z_r$ defines the value of the Walsh function $W_k(x) $ on each interval $x \in ((\ell-1)/n, \ell/n)$ for each $\ell=1,\ldots,n$ where $n = 2^r$.  The disadvantage is that there is no efficient direct numerical construction.  Thus the Walsh matrices with the sequency ordering are normally constructed by permutation of the natural ordering.  We will always use the natural order in this paper.

\subsection{Contribution}
\label{con}

In this paper we discuss discrepancy results for a special class of Parseval frames defined by Walsh matrices.  In particular we show that if $m,n \in {\mathbb N}$ with $m \leq n = 2^r$ for some $r \in {\mathbb N}$ then there is a Parseval frame defined by a pre-frame matrix operator $V = [\bfv_1,\ldots,\bfv_n] \in {\mathbb C}^{m \times n}$ where $\bfv_j = [v_{ij}] \in {\mathbb C}^m$ with $v_{ij} = \pm 1/\sqrt{n}$ for each $i=1,\ldots,m$ and $j=1,\ldots,n$ and $\| \bfv_j \| = \sqrt{m/n}$ for each $j=1,\ldots,n$.  We show that for $m \leq n/2$ these frames can be split into two identical tight frames with frame constant $c =1/2$.  For $n/2 < m \leq n$ we show that the frames can no longer be evenly split but we find an explicit expression for the discrepancy in a best possible split.  Because the vectors in our frames are all the same length we have not imposed any direct condition that forces them to be small.  Hence our results are not directly comparable to those in MSSDT.  Of course the frame vectors are small if $\sqrt{m/n}$ is small.   We also show that all Parseval frames in ${\mathbb C}^m$ constructed from $n = 2^r$ vectors of equal length can be transformed to a corresponding Walsh frame and we ponder the implications of this correspondence in regard to splitting of the associated quadratic forms.  

\section{The main results}
\label{mr}

Let $m,n \in {\mathbb N}$ with $m \leq n$.  We would like to construct a Parseval frame defined by a pre-frame matrix operator $V = [ v_{ij}] = [\bfv_1,\ldots,\bfv_n] \in {\mathbb C}^{m \times n}$ where the frame vectors $\bfv_j \in {\mathbb C}^m$ all have the same length.  In order to construct the simplest possible frame we will insist that $v_{ij} = \pm\, 1 /\sqrt{n}$ for all $i=1,\ldots,m$ and all $j = 1,\ldots,n$.  Thus $\| \bfv_j\| = \sqrt{m/n}$ for all $j=1,\ldots,n$.  To facilitate splitting the frame into two potentially equal parts we will restrict our attention to frames with $n = 2^r$ vectors for some $r \in {\mathbb N}$.  We will show that the normalized Walsh matrices provide the ideal building blocks for our proposed frames.  We discuss splitting of Parseval frames defined by Walsh matrices and find explicit expressions for the minimal discrepancy.

\subsection{Parseval frames defined by Walsh matrices}
\label{ntfwm}

Define $n = 2^r$ for some $r \in {\mathbb N} + 1$ and suppose $ m \in {\mathbb N}$ with $m < n$.  Thus we exclude the case $m=n$.  Let $Y = [\bfy_1,\ldots,\bfy_n] \in {\mathbb C}^{n \times n}$ be the corresponding Walsh matrix.  We have
$$
Y^*Y = YY^* = nI_n
$$
where $I_n \in {\mathbb C}^{n \times n}$ is the unit matrix and so $G = Y/\sqrt{n}$ is a unitary matrix.  Write $G = [\bfg_1,\ldots,\bfg_m \mid \bfg_{m+1},\ldots,\bfg_n]$ and define $W \in {\mathbb C}^{n \times m}$ by setting $W = [\bfg_1,\ldots,\bfg_m]$.  Let $V = W^* \in {\mathbb C}^{m \times n}$.  Since the columns of $W$ are a subset of the columns of $G$ they form an orthonormal set in ${\mathbb C}^n$.  Therefore $W^*W = I_m \iff VV^* = I_m$.  As usual we write $V = [\bfv_1,\ldots,\bfv_n]$ where $\bfv_j \in {\mathbb C}^m$ for all $j=1,\ldots,n$.  The column vectors $\{\bfv_1,\ldots,\bfv_n\}$ form a Parseval frame in ${\mathbb C}^m$ and since $v_{ij} = \pm1/\sqrt{n}$ for each $i=1,\ldots,m$ and $j=1,\ldots,n$ it follows that $\|\bfv_j\| = \sqrt{m/n}$ for all $j=1,\ldots,n$.   Thus the frame vectors are all the same size.  The Parseval frame defined by the columns $V = [\bfv_1,\ldots,\bfv_n] \in {\mathbb R}^{m \times n}$ will be called a {\em Walsh frame}.

We wish to consider what happens when we try to split a Walsh frame into two equal parts.  We begin with a simple example.

\begin{example}
\label{ex1}

{\rm Let $m=3$ and $n = 8$.  Use the first three rows of the Walsh matrix $Y_3$ to define
$$
V = \frac{1}{2\sqrt{2}} \left[ \begin{array}{rrrr|rrrr}
1 & 1 & 1 & 1 & 1 & 1 & 1 & 1 \\
1 & -1 & -1 & 1 & 1 & -1 & -1 & 1 \\
1 & 1 & -1 & -1 & 1 & 1 & -1 & -1 \end{array} \right] = \left[ \begin{array}{c|c}
V_a & V_b \end{array} \right]
$$
so that $VV^* = I_3$.  The Parseval frame can be split into two identical frames as shown above with $V_a = V_b$ and $V_a{V_a}^* = V_b{V_b}^*= I_3/2$.  Now renormalize and define
$$
V = \frac{1}{2} \left[ \begin{array}{rr|rr}
1 & 1 &1 & 1 \\
1 & -1 & -1 & 1 \\
1 & 1 & -1 & -1 \end{array} \right] = \left[ \begin{array}{c|c}
V_a & V_b \end{array} \right]
$$
so that $VV^* = I_3$.  If we split the new Parseval frame into two parts $V_a$ and $V_b$ as shown above then the two parts are no longer identical.  In fact a little thought will show that no even split is possible.  For the proposed split we have
$$
\|V_a{V_a}^* - I_3/2 \|_2 = \|V_b{V_b}^* - I_3/2 \|_2 = 1/2
$$
which is the best possible.  For the corresponding quadratic forms we have 
\begin{eqnarray*}
s(\bfx) & = & \bfx^*VV^*\bfx \\
& = & x_1^2 + x_2^2 + x_3^2 \\
& = & \left( x_1^2/2 + x_2^2/2 + x_3^2/2 + x_1x_3 \right) + \left( x_1^2/2 + x_2^2/2 + x_3^2/2 - x_1x_3 \right) \\
& = & \bfx^*V_a{V_a}^*\bfx + \bfx^*V_b{V_b}^*\bfx \\
& = & s_a(\bfx) + s_b(\bfx).
\end{eqnarray*}
Considered separately the sets of vectors defined by the columns of $V_a$ and $V_b$ no longer span ${\mathbb C}^3$.  Thus neither $V_a$ nor $V_b$ defines a frame for ${\mathbb C}^3$.}  $\hfill \Box$
\end{example}

Let $m \in {\mathbb N}$ with $2^{s-1} < m \leq 2^s$ for some $s \in {\mathbb N}$ and let $n = 2^r$ for some $r \in {\mathbb N}+s$.  Thus $m < 2^{s+1} \leq n$.   Consider a Walsh frame defined by the first $m$ rows of the normalized Walsh matrix $F = Y/\sqrt{n}$ where $Y = Y_r$ is the corresponding Walsh matrix.  The above example suggests that we can split this Parseval frame into two identical tight sub-frames each having $2^{r-1}$ elements.  Indeed the example suggests that we can split the tight frame into identical tight sub-frames $k$ times where $k = r-s$.  We have the following elementary result.

\begin{theorem}[WF1]
\label{wft1}
Let $n=2^r$ for some $r \in {\mathbb N}$ and suppose $m \in {\mathbb N}$ with $m \leq n/2 = 2^{r-1}$.  Let $W = [ \bfy_1,\ldots,\bfy_m ]/\sqrt{n} \in {\mathbb C}^{n \times m}$ be defined by the first $m$ columns of the Walsh matrix $Y_r \in {\mathbb C}^{n \times n}$ and let $V = [\bfv_1,\ldots,\bfv_n] = W^* \in {\mathbb C}^{m \times n}$.  Then $VV^* = I_m$ and the Parseval frame for ${\mathbb C}^m$ defined by the columns of the matrix $V$ can be split into two identical tight frames for ${\mathbb C}^m$ defined by the columns of the matrices $V_1 = [\bfv_1,\ldots,\bfv_{n/2}]$ and $V_2 = [\bfv_{n/2+1},\ldots,\bfv_n]$ with $V_1{V_1}^* = V_2{V_2}^* = I_m/2$. $\hfill \Box$
\end{theorem} 

\noindent{\bf Proof} \quad  It follows from the recursive definition of the Walsh matrices 
$$
Y_r = \left[ \begin{array}{rr}
Y_{r-1} & Y_{r-1} \\
Y_{r-1} & -Y_{r-1} \end{array} \right]
$$
that $V_1, V_2 \in {\mathbb C}^{m \times n/2}$ are sub-matrices of ${Y_{r-1}}^*/\sqrt{n} = Y_{r-1}/\sqrt{n}$ consisting in each case of the first $m$ rows.  Hence they are identical.  Each of the matrices $V_1, V_2$ has $m \leq n/2$ mutually orthogonal rows and each row has length $1/\sqrt{2}$.  $\hfill \Box$ \bigskip

Although {\em redundancy} and the additional associated {\em flexibility} are useful ingredients in the use of frames \cite{cas2,chr1} the redundancy in the Parseval frames defined by Walsh matrices is simply repetition of individual vectors.  If $m = 2^k < 2^r = n$ then each vector is repeated $2^{k - r}$ times.  Thus an equal split is obvious.  We will now restrict our attention to reduced Walsh frames with $n = 2^r$ and $n/2 < m \leq n$.   If $m = n/2 +s$ then $s$ vectors from the $n/2$ individual vectors in the frame are repeated.  A reduced Walsh frame cannot be evenly split.  Define $Y \in {\mathbb C}^{n \times m}$ by deleting $n/2 - s$ arbitrarily-selected columns from the right-hand half of the Walsh matrix $Y_r$.  Thus we have
$$
Y = \left[ \begin{array}{rrr|rrr}
\bfy_1 & \cdots & \bfy_{n/2} & \bfy_{k(1)} & \cdots & \bfy_{k(s)} \\ \hline
\bfy_1 & \cdots & \bfy_{n/2} & - \bfy_{k(1)} & \cdots & -\bfy_{k(s)} \end{array} \right] 
$$
where $Y_{r-1} = [\bfy_1,\ldots,\bfy_{n/2}] \in {\mathbb C}^{(n/2) \times (n/2)}$ is the Walsh matrix of order $r-1$ and $1 \leq k(1) < \cdots < k(s) \leq n/2$ is an arbitrarily-selected subset of size $s$ from $\{1,\ldots,n/2\}$.  If we define $W = Y/ \sqrt{n}$ and $V = W^*$ then we have
$$
VV^* = W^*W = I.
$$
Hence the columns of $V$ form a Parseval frame for ${\mathbb C}^m$ with $\|\bfv_j\| = \sqrt{m/n}$ for all $j=1,\ldots,n$.   We wish to split the frame as evenly as possible.  Let
\begin{eqnarray*}
Y_a & = & \left[ \begin{array}{rrr|rrr}
\bfy_1 & \cdots & \bfy_{n/2} & \bfy_{k(1)} & \cdots & \bfy_{k(s)}  \end{array} \right] \in {\mathbb C}^{(n/2) \times m} \\
Y_b & = & \left[ \begin{array}{rrr|rrr}
\bfy_1 & \cdots & \bfy_{n/2} & - \bfy_{k(1)} & \cdots & -\bfy_{k(s)} \end{array} \right] \in {\mathbb C}^{(n/2) \times m}
\end{eqnarray*}
and define $W_a = Y_a / \sqrt{n}$, $W_b = Y_b / \sqrt{n}$, $V_a = {W_a}^*$ and $V_b = {W_b}^*$.  We know that
$$
\bfy_j^* \bfy_k = \frac{n}{2}\, \delta_{jk} = \left\{ \begin{array}{ll}
\displaystyle \frac{n}{2} &\mbox{if}\ k = j \\
& \\
0 &\mbox{otherwise}. \end{array} \right.
$$
Therefore
\begin{eqnarray*}
{Y_a}^*Y_a & = & \left[ \begin{array}{c}
\left[ \begin{array}{c} \bfy_j^* \end{array} \right] \\ 
\vspace{-0.375cm} \\ \hline
\vspace{-0.375cm} \\
\left[ \begin{array}{c} \bfy_{k(p)}^*\end{array} \right] \end{array} \right] \left[ \rule{0cm}{0.4cm} \begin{array}{c|c}
\left[ \begin{array}{c} \bfy_k \end{array} \right] & \left[ \begin{array}{c} \bfy_{k(q)} \end{array} \right] \end{array} \right] \\
& = & \left[ \begin{array}{c|c}
\left[ \begin{array}{c} \bfy_j^* \bfy_k \end{array} \right] & \left[ \begin{array}{c} \bfy_j^* \bfy_{k(q)} \end{array} \right] \\ 
\vspace{-0.375cm} \\ \hline
\vspace{-0.375cm} & \\
\left[ \begin{array}{c} \bfy_{k(p)}^* \bfy_k \end{array} \right] & \left[ \begin{array}{c} \bfy_{k(p)} \bfy_{k(q)} \end{array} \right] \end{array} \right] \\
& = & \frac{n}{2} \left[ \begin{array}{c|c}
\left[ \begin{array}{c} \delta_{jk} \end{array} \right] & \left[ \begin{array}{c} \delta_{j k(q)} \end{array} \right] \\ 
\vspace{-0.375cm} \\ \hline
\vspace{-0.375cm} & \\ 
\left[ \begin{array}{c} \delta_{k(p)k} \end{array} \right] & \left[ \begin{array}{c} \delta_{k(p)k(q)} \end{array} \right] \end{array} \right] \\
& = & \frac{n}{2} \left[ \begin{array}{c|c}
I_{n/2} & \Delta \\
\vspace{-0.375cm} \\ \hline
\Delta^* & I_s \end{array} \right] 
\end{eqnarray*}
where $\Delta = [ \bfe_{k(1)}, \ldots, \bfe_{k(s)} ] \in {\mathbb C}^{(n/2) \times s}$.  A similar argument shows that
$$
{Y_b}^*Y_b = \frac{n}{2} \left[ \begin{array}{c|c}
I_{n/2} & - \Delta \\
\vspace{-0.35cm} \\ \hline
- \Delta^* & I_s \end{array} \right]. 
$$
Thus
$$
V_a{V_a}^* = \frac{1}{2} \left[ \begin{array}{c|c}
I_{n/2} & \Delta \\ \hline
\Delta^* & I_s \end{array} \right] \quad \mbox{and} \quad V_b{V_b}^* = \frac{1}{2} \left[ \begin{array}{c|c}
I_{n/2} & - \Delta \\ \hline
- \Delta^* & I_s \end{array} \right].
$$
Our chosen split is not even but is nevertheless the best possible.  What is the discrepancy in this case?  Let
$$
E_a = V_a{V_a}^* - I_m/2 = \frac{1}{2} \left[ \begin{array}{cc}
0 & \Delta \\
\Delta^* & 0 \end{array} \right]
$$
and
$$
E_b = V_b{V_b}^* - I_m/2 = \frac{1}{2} \left[ \begin{array}{cc}
0 & - \Delta \\
- \Delta^* & 0 \end{array} \right].
$$
We have $\Delta^* \Delta = I_s$ and $\Delta \Delta^* = U = [u_{pq}] \in {\mathbb C}^{(n/2) \times (n/2)}$  where $u_{pq} = 1$ if $p=q=k(r)$ for $r=1,\ldots,s$ and $u_{pq} = 0$ otherwise.  Therefore
$$
E_a^2 = E_b^2 = \frac{1}{4} \left[ \begin{array}{cc}
U & 0 \\
0 & I_s \end{array} \right]
$$
and hence the eigenvalues of $E_a^2$ and $E_b^2$ are $\lambda_1 = 0$ with multiplicity $n/2 - s$ corresponding to the zero rows and columns and $\lambda_2 = 1/4$ with multiplicity $2s$ corresponding to unit rows and columns.  Therefore $\|E_a^2\|_2 = \|E_b^2\|_2 = 1/4$ and hence $\|E_a\|_2 = \|E_b\|_2 = 1/2$.

Alternatively we can show that the eigenvalues of the real symmetric matrix $V_a{V_a}^*$ are $\lambda = 0, 1/2, 1$.   We have
$$
\lambda I_m - V_a{V_a}^* = \left[ \begin{array}{cc}
(\lambda - 1/2)I_{n/2}  & - \Delta/2 \\
- \Delta^*/2 & (\lambda - 1/2) I_s \end{array} \right].
$$
We can use left multiplication by elementary matrices to perform elementary row operations on the matrix $\lambda I_m - V_a{V_a}^*$ and thereby show that $\lambda = 0$ and $\lambda = 1$ are each eigenvalues of multiplicity $s$.  For $\lambda = 0$ we can see that
$$
\left[ \begin{array}{cc}
I_{n/2} & 0 \\
- \Delta^* & I_s \end{array} \right] \left[ \begin{array}{cc}
- I_{n/2} & - \Delta \\
- \Delta^* & - I_s \end{array} \right] = \left[ \begin{array}{cc}
- I_{n/2} & - \Delta \\
0 & \Delta^* \Delta - I_s \end{array} \right] = \left[ \begin{array}{cc}
- I_{n/2} & - \Delta \\
0 & 0 \end{array} \right]
$$
and for $\lambda = 1$ we have    
$$
\left[ \begin{array}{cc}
I_{n/2} & 0 \\
\Delta^* & I_s \end{array} \right] \left[ \begin{array}{cc}
I_{n/2} & - \Delta \\
- \Delta^* & I_s \end{array} \right] = \left[ \begin{array}{cc}
I_{n/2} & - \Delta \\
0 & - \Delta^* \Delta + I_s \end{array} \right] = \left[ \begin{array}{cc}
I_{n/2} & - \Delta \\
0 & 0 \end{array} \right].
$$
In each case the reduced matrix has rank $n/2$ and hence the eigenvalue has multiplicity $m - n/2 = s$.  For $\lambda = 1/2$ the matrix
$$
 I_m/2 - V_a{V_a}^* = \left[ \begin{array}{cc}
0  & - \Delta/2 \\
- \Delta^*/2 & 0 \end{array} \right]
$$
has rank $2s$ and so $\lambda = 1/2$ is an eigenvalue of multiplicity $m - 2s$.  It follows that $\| V_a{V_a}^*\|_2 = \| {V_a}^*V_a \|_2 = 1$.  A similar argument shows us that $\| V_b{V_b}^*\|_2 = \| {V_b}^*V_b \|_2 = 1$.  If we use the Frobenius norm then
$$
\| V_a{V_a}^* - I_m/2 \|_F = \| V_b{V_b}^* - I_m/2 \|_F = \sqrt{ 2 \| \Delta/2\|_F^2} = \sqrt{s/2} = \sqrt{(m - n/2)/2}.
$$
This leads us to the following result.

\begin{theorem}[WF2]
\label{wft2}
Let $n=2^r$ for some $r \in {\mathbb N}$ and suppose $m \in {\mathbb N}$ with $2^{r-1} = n/2 \leq m < n = 2^r$.  Let $Y_{r-1} = [\bfy_1,\ldots,\bfy_{n/2}] \in {\mathbb C}^{(n/2) \times (n/2)}$ be the Walsh matrix of order $r-1$ and define  
$$
Y = \left[ \begin{array}{rrr|rrr}
\bfy_1 & \cdots & \bfy_{n/2} & \bfy_{k(1)} & \cdots & \bfy_{k(s)} \\ \hline
\bfy_1 & \cdots & \bfy_{n/2} & - \bfy_{k(1)} & \cdots & -\bfy_{k(s)} \end{array} \right] = \left[ \begin{array}{cc}
Y_a \\ \hline
Y_b \end{array} \right] \in {\mathbb C}^{n \times m}
$$
where $s = m - n/2$ and $1 \leq k(1) < \cdots < k(s) \leq n/2$ is an arbitrarily selected subset of $\{1,\ldots,n/2\}$.  Let $V = Y^*/\sqrt{n}$, $V_a = {Y_a}^*/\sqrt{n}$ and $V_b = {Y_b}^*/\sqrt{n}$.  Then $VV^* = I$ and the split defined by $V = [V_a \mid V_b]$ gives
$$
V_a{V_a}^* - I/2  = \left[ \begin{array}{c|c}
0 & \Delta/2 \\ \hline
\Delta^*/2 & 0 \end{array} \right] \quad \mbox{and} \quad V_b{V_b}^* - I/2 = \left[ \begin{array}{c|c}
0 & - \Delta/2 \\ \hline
- \Delta^*/2 & 0 \end{array} \right]
$$
where $\Delta = [\bfe_{k(1)} \cdots \bfe_{k(s)}] \in {\mathbb C}^{(n/2) \times s}$.  The error
$$
\|V_a{V_a}^* - I/2\|_2 = \|V_b{V_b}^* - I/2\|_2 = 1/2
$$
is the best possible.  We also have $\| V_a{V_a}^*\|_2 = \|V_b{V_b}^*\|_2 = 1$.  If we use the Frobenius norm then
$$
\|V_a{V_a}^* - I/2\|_F = \|V_b{V_b}^* - I/2\|_F = \sqrt{s/2} = \sqrt{(m - n/2)/2}
$$
and $\| V_a{V_a}^*\|_F = \|V_b{V_b}^*\|_F = \sqrt{(m + 2s)/4} = \sqrt{(3m-n)/4}$.  $\hfill \Box$
\end{theorem}

We illustrate our results with an example.

\begin{example}
\label{ex2}

{\rm In the case where $m=6$ and $r = 3$ we have $n = 2^3 = 8$ with
$$
Y_2 = \left[ \begin{array}{cccc}
\bfy_1 & \bfy_2 & \bfy_3 & \bfy_4 \end{array} \right] = \left[ \begin{array}{rrrr}
1 & 1 & 1 & 1  \\
1 & -1 & 1 & -1 \\
1 & 1 & -1 & -1 \\
1 & -1 & -1 & 1 \end{array} \right].
$$
If we choose $\bfk = [ 1,\  3 ]$ then we define
\begin{eqnarray*}
Y & = & \left[ \begin{array}{rrrr|rr}
\bfy_1 & \bfy_2 & \bfy_3 & \bfy_4 & \bfy_1 & \bfy_3 \\ \hline
\bfy_1 & \bfy_2 & \bfy_3 & \bfy_4 & -\bfy_1 & -\bfy_3 \end{array} \right] \\
& = & \left[ \begin{array}{rrrr|rr}
1 & 1 & 1 & 1 & 1 & 1 \\
1 & -1 & 1 & -1 & 1 & 1 \\
1 & 1 & -1 & -1 & 1 & -1 \\
1 & -1 & -1 & 1 & 1 & -1 \\ \hline
1 & 1 & 1 & 1 & -1 & -1 \\
1 & -1 & 1 & -1 & -1 & -1 \\
1 & 1 & -1 & -1 & -1 & 1 \\
1 & -1 & -1 & 1 & -1 & 1 \end{array} \right]
\end{eqnarray*}
and $V = Y^*/(2 \sqrt{2})$.  Now define 
\begin{eqnarray*}
Y_a & = & \left[ \begin{array}{rrrr|rr}
\bfy_1 & \bfy_2 & \bfy_3 & \bfy_4 & \bfy_1 & \bfy_3 \end{array} \right] \\
& & \\
& = & \left[ \begin{array}{rrrr|rr}
1 & 1 & 1 & 1 & 1 & 1 \\
1 & -1 & 1 & -1 & 1 & 1 \\
1 & 1 & -1 & -1 & 1 & -1 \\
1 & -1 & -1 & 1 & 1 & -1 \end{array} \right]
\end{eqnarray*}
and
\begin{eqnarray*}
Y_b & = & \left[ \begin{array}{rrrr|rr}
\bfy_1 & \bfy_2 & \bfy_3 & \bfy_4 & - \bfy_1 & - \bfy_3 \end{array} \right]
\\
& & \\
& = & \left[ \begin{array}{rrrr|rr}
1 & 1 & 1 & 1 & -1 & -1 \\
1 & -1 & 1 & -1 & -1 & -1 \\
1 & 1 & -1 & -1 & -1 & 1 \\
1 & -1 & -1 & 1 & -1 & 1 \end{array} \right]
\end{eqnarray*}
and let $V_a = {Y_a}^*/(2\sqrt{2})$ and $V_b = {Y_b}^*/(2\sqrt{2})$.  We have 
$$
V_a{V_a}^* =  \left[ \begin{array}{rrrr|rr}
\rule{0cm}{0.4cm} \frac{1}{2} & 0 & 0 & 0 & \frac{1}{2} & 0 \\
\rule{0cm}{0.4cm} 0 & \frac{1}{2} & 0 & 0 & 0 & 0 \\
\rule{0cm}{0.4cm} 0 & 0 & \frac{1}{2} & 0 & 0 & \frac{1}{2} \\
\rule{0cm}{0.4cm} 0 & 0 & 0 & \frac{1}{2} & 0 & 0 \\ 
\vspace{-0.375cm} &&&& \\ \hline
\rule{0cm}{0.4cm} \frac{1}{2} & 0 & 0 & 0 & \frac{1}{2} & 0 \\
\rule{0cm}{0.4cm} 0 & 0 & \frac{1}{2} & 0 & 0 & \frac{1}{2} \end{array} \right] = \left[ \begin{array}{c|c}
I_4/2 & \Delta/2 \\ \hline
\Delta^*/2 & I_2/2 \end{array} \right],
$$
and
$$
V_b{V_b}^* =  \left[ \begin{array}{rrrr|rr}
\rule{0cm}{0.4cm} \frac{1}{2} & 0 & 0 & 0 & -\frac{1}{2} & 0 \\
\rule{0cm}{0.4cm} 0 & \frac{1}{2} & 0 & 0 & 0 & 0 \\
\rule{0cm}{0.4cm} 0 & 0 & \frac{1}{2} & 0 & 0 & -\frac{1}{2} \\
\rule{0cm}{0.4cm} 0 & 0 & 0 & \frac{1}{2} & 0 & 0 \\
\vspace{-0.375cm} &&&& \\ \hline
\rule{0cm}{0.4cm} -\frac{1}{2} & 0 & 0 & 0 & \frac{1}{2} & 0 \\
\rule{0cm}{0.4cm} 0 & 0 & -\frac{1}{2} & 0 & 0 & \frac{1}{2} \end{array} \right] = \left[ \begin{array}{c|c}
I_4/2 & - \Delta/2 \\ \hline
- \Delta^*/2 & I_2/2 \end{array} \right],
$$
where
$$
\Delta = \left[ \begin{array}{rr}
1 & 0 \\
0 & 0 \\
0 & 1 \\
0 & 0 \end{array} \right].
$$
Note that $\Delta^*\Delta = I_2$ and
$$
U = \Delta \Delta^* = \left[ \begin{array}{cccc}
1 & 0 & 0 & 0 \\
0 & 0 & 0 & 0 \\
0 & 0 & 1 & 0 \\
0 & 0 & 0 & 0 \end{array}\right].
$$
We have
$$
\| V_a{V_a}^* - I_4/2\|_2 = \| V_b{V_b}^* - I_4/2\|_2 = 1/2
$$
and $\| V_a{V_a}^* \| = \| V_b{V_b}^*\| = 1$.  If we return to the idea that each matrix $S = VV^*$ is a sum of elementary rank $1$ matrices then we have $S = \sum_{j=1}^6 P_j$ where $P_j = \bfv_j {\bfv_j}^* = \bfy_j{\bfy_j}^*/8$ for each $j=1,\ldots,6$.  If we use the Frobenius norm then $\|P_j\|_F^2 = 1/4$ and
$$
\| V_a{V_a}^* - I_4/2\|_F = \| V_b{V_b}^* - I_4/2\|_F = 1.
$$
These calculations agree with the general results stated in WF2.}  $\hfill \Box$ 
\end{example}

\section{Parseval frames defined by vectors of equal length}
\label{ntfelv}

The quadratic forms defined by Walsh frames can be split exactly if $m \leq n/2 = 2^{r-1}$.  If $n/2 < m < n = 2^r$ the quadratic forms can no longer be evenly split but there is an explicit description for the minimal discrepancy.  An interesting question is whether these results are completely specific to Walsh frames or whether similar results apply to quadratic forms defined by other Parseval frames constructed from vectors of equal length.  We begin by stating a well-known lemma.

\begin{lemma}
\label{frameconstant}
Let $m \in {\mathbb N}$ and $n = 2^r$ for some $r \in {\mathbb N}$.  If $m \leq n$ and $V = [\bfv_1,\ldots,\bfv_n] \in {\mathbb C}^{m \times n}$ defines a Parseval frame for ${\mathbb C}^m$ with $S = VV^* = I_m$ and $\|\bfv_j\|^2 = \alpha$ for each $j=1,\ldots,m$ then $\alpha = m/n$. $\hfill \Box$
\end{lemma}

Let $m \in {\mathbb N}$ with $m \leq n = 2^r$ for some $r \in {\mathbb N}$.  Suppose that $V = [\bfv_1,\ldots,\bfv_n] \in {\mathbb C}^{m \times n}$ defines a Parseval frame with $\|\bfv_j\| = \sqrt{m/n}$ for each $j=1,\ldots,m$.  Define $W = [\bfw_1,\ldots,\bfw_m] \in {\mathbb C}^{n \times m}$ by setting $W = V^*$.  We have $W^*W = VV^* = I_m$ and so $\{\bfw_1,\ldots,\bfw_m\} \in {\mathbb C}^n$ is an orthonormal set.   Let us extend this set to an orthonormal basis $\{\bfw_1,\ldots,\bfw_n\} \in {\mathbb C}^n$.  Define the orthogonal matrix
$$
H = [H_1 \mid H_2 ] = [ \bfw_1,\ldots,\bfw_m \mid \bfw_{m+1},\ldots, \bfw_n] \in {\mathbb C}^{n \times n}
$$
where $H_1 = W \in {\mathbb C}^{n \times m}$ and $H_2 \in {\mathbb C}^{n \times (n-m)}$.  Define $G = H^* \in {\mathbb C}^{n \times n}$.  We can write
$$
G = \left[ \begin{array}{c}
G_1 \\ \hline
G_2 \end{array} \right] = \left[ \begin{array}{cccc}
\bfv_1 & \bfv_2 & \cdots & \bfv_n \\ \hline
\bfr_1 & \bfr_2 & \cdots & \bfr_n \end{array} \right]
$$
where $G_1 = V \in {\mathbb C}^{m \times n}$ and $G_2 = R = [\bfr_1,\ldots,\bfr_n] \in {\mathbb C}^{(n-m) \times n}$.  If we define
$$
\bfg_j = \left[ \begin{array}{c}
\bfv_j \\
\bfr_j \end{array} \right]
$$
for each $j=1,\ldots,n$ then we can write $G = [\bfg_1,\ldots,\bfg_n]$.  The matrix $G$ defines an orthonormal basis $\{\bfg_1,\ldots,\bfg_n\}$ for ${\mathbb C}^n$.  The set $\{\bfg_1,\ldots,\bfg_n\}$ also defines an embedded Parseval frame for the $m$-dimensional subspace $S_m$ of ${\mathbb C}^n$ spanned by all vectors of the form
$$
\bfz = \left[ \begin{array}{c}
\bfx \\ 
\bfzero \end{array} \right]
$$
where $\bfx \in {\mathbb C}^m$.   Let $Y = [\bfy_1,\ldots,\bfy_n] \in {\mathbb C}^{n \times n}$ be the Walsh matrix of order $r$ and define $F = [\bff_1,\ldots,\bff_n] \in {\mathbb C}^{n \times n}$ by setting $F = Y/\sqrt{n} \in {\mathbb C}^{n \times n}$.  We have $\bff_j = \bfy_j/\sqrt{n}$ for all $j=1,\ldots,n$.  Note that the normalized Walsh matrix $F$ is real symmetric and orthogonal.  Define an orthogonal matrix $P \in {\mathbb C}^{n \times n}$ by setting $P = F H$.   Therefore $PG = F$ and hence $P \bfg_j = \bff_j$ for all $j=1,\ldots,n$.  We will use the orthogonal matrix $P$ to change the coordinate representation for the embedded Parseval frame defined by $G$ into a representation defined by $F$.   Thus the embedded frame now looks like a Walsh frame.  

We saw earlier that we can represent vectors in $S_m$ using the embedded frame with
\begin{equation}
\label{genrep2}
\left[ \begin{array}{c}
\bfx \\ 
\bfzero \end{array} \right] = \sum_{j=1}^n ( {\bfg_j}^* \bfz) \bfg_j = \sum_{j=1}^n ( {\bfv_j}^* \bfx) \left[ \begin{array}{c}
\bfv_j \\ 
\bfr_j \end{array} \right]
\end{equation}
for all $\bfx \in {\mathbb C}^m$.  To see this representation in the new coordinates we simply multiply both sides of (\ref{genrep2}) by $P$.  Thus we have
\begin{equation}
\label{genrep3}
\bfy = P \left[ \begin{array}{c}
\bfx \\
\bfzero \end{array} \right] = \sum_{j=1}^n ( {\bfg_j}^* \bfz)  P \bfg_j = \sum_{j=1}^n ( {\bfv_j}^* \bfx) \bff_j
\end{equation}
where we have the same coefficients yet again.  Since $VV^* = I_m$ the quadratic form for the original frame is simply $$
s(\bfx) = \bfx^* \bfx = \sum_{j=1}^m |{\bfv_j}^* \bfx |^2 = \sum_{j=1}^n s_j(\bfx)
$$
where the $s_j(\bfx)$ are rank $1$ quadratic forms for all $j=1,\ldots,n$.  Since $GG^* = I_n$ the corresponding embedded quadratic form is given by
$$
t(\bfz) = \bfz^* \bfz = \sum_{j=1}^n | {\bfg_j}^* \bfz |^2 = \sum_{j=1}^n t_j(\bfz)
$$
where
$$
\bfz = \left[ \begin{array}{c}
\bfx \\
\bfr \end{array} \right] \in {\mathbb C}^n
$$
and where the $t_j(\bfz)$ are rank $1$ quadratic forms.  For quadratic forms on the subspace $S_m$ we have
$$
t \left( \left[ \begin{array}{c}
\bfx \\ 
\bfzero \end{array} \right] \right) = \bfx^* \bfx + \bfzero^* \bfzero = \sum_{j=1}^n | {\bfv_j}^* \bfx + {\bfr_j }^* \bfzero |^2 = \sum_{j=1}^m | {\bfv_j}^* \bfx|^2 = \sum_{j=1}^m s_j(\bfx) = s(\bfx).
$$
In the new coordinates $\bfy = P^* \bfz \iff P \bfy = \bfz$ we can use (\ref{genrep3}) when $\bfr = \bfzero$ to see that
$$
q(\bfy) = \bfy^* \bfy = \sum_{j=1}^n \sum_{k=1}^n \bff_j^*(\bfx^* \bfv_j)( {\bfv_k}^*\bfx) \bff_k = \sum_{j=1}^n |{\bfv_j}^* \bfx |^2 = s(\bfx).
$$
If we write 
$$
P = \left[ \begin{array}{c}
P_1 \\
P_2 \end{array} \right]
$$
where $P_1 \in {\mathbb C}^{m \times n}$ and $P_2 \in {\mathbb C}^{(n-m) \times n}$ then the subspace $S_m$ defined by $\bfr = \bfzero$ is defined in the new coordinates by $P_2 \bfy = \bfzero$.  Since $P$ is invertible we must have $\mbox{rank}(P_2) = n-m$.  Thus we can use elementary row operations to eliminate $(n-m)$ variables and hence express $q(\bfy)$ as a sum of rank $1$ quadratic forms in $m$ variables $y_{k(1)},y_{k(2)},\ldots,y_{k(m)}$ where $1 \leq k(1) < k(2) < \cdots < k(m) \leq n$. 

Although we have assumed $n=2^r$ throughout we shall see in the following example that this assumption is basically just a matter of convenience.  If $V = [\bfv_1,\ldots,\bfv_k] \in {\mathbb C}^{m \times k}$ form a Parseval frame in ${\mathbb C}^k$ where $2^{r-1} = n/2 < k < n = 2^r$ and if we define $W = [\bfw_1,\ldots,\bfw_m] \in {\mathbb C}^{k \times m}$ by setting $W = V^*$ then the columns $\{\bfw_1,\ldots,\bfw_m\}$ form an orthonormal set in ${\mathbb C}^k$.  We can embed these vectors in ${\mathbb C}^n$ by defining
$$
\bfh_j = \left[ \begin{array}{c}
\bfw_j \\
\bfzero \end{array} \right] \in {\mathbb C}^n
$$
for each $j=1,\ldots,m$.  Now $\{\bfh_1,\ldots,\bfh_m\}$ forms an orthonormal set in ${\mathbb C}^n$ which we can easily extend to an orthonormal basis $\{\bfh_1,\ldots,\bfh_n\} \in {\mathbb C}^n$.  Define orthogonal matrices $H = [\bfh_1,\ldots,\bfh_n] \in {\mathbb C}^{n \times n}$ and $G = [\bfg_1,\ldots,\bfg_n] \in {\mathbb C}^{n \times n}$ by setting $G = H^*$.  The vectors $\{\bfg_1,\ldots,\bfg_n\} \in {\mathbb C}^n$ now form a Parseval frame for the $m$-dimensional subspace $S_m$ defined by vectors in the form
$$
\bfz = \left[ \begin{array}{c}
\bfx \\ 
\bfzero \end{array} \right] \in {\mathbb C}^n
$$
for all $\bfx \in {\mathbb C}^m$.  The frame defined by $G$ for the $m$-dimensional subspace is simply the frame defined by $V$ embedded into ${\mathbb C}^n$.  We can write $G \cong V$.
  
We have argued that from within ${\mathbb C}^n$ all orthonormal bases look the same and that coordinate representation is essentially a matter of choice.  Thus it is always possible to use the columns of a normalized Walsh matrix $F$ to represent the vectors of a Parseval frame defined by vectors of equal length.  We illustrate our remarks by considering a particular example. 

\begin{example}
\label{ex3}

{\rm Suppose we wish to find three vectors of equal length that form a tight frame in ${\mathbb C}^2$.  If we define
$$
\bfv_1 = \left[ \begin{array}{cc}
a \\
\sqrt{t^2 - a^2} \end{array} \right], \quad \bfv_2 = \left[ \begin{array}{cc}
b \\
\sqrt{t^2 - b^2} \end{array} \right],  \quad \bfv_3 = \left[ \begin{array}{cc}
c \\
\sqrt{t^2 - c^2} \end{array} \right]
$$
then $\| \bfv_j\| = |t|$ for each $j=1,2,3$.  The condition for a Parseval frame is that
$$
\bfw_1 = \left[ \begin{array}{c}
a \\
b \\
c \end{array} \right] \quad \mbox{and} \quad \bfw_2 = \left[ \begin{array}{c}
\sqrt{t^2 - a^2} \\
\sqrt{t^2 - b^2} \\
\sqrt{t^2 - c^2} \end{array} \right] 
$$
form an orthonormal set.  Thus we require $a^2 + b^2 + c^2 = 1$, $(t^2 - a^2) + (t^2 - b^2) + (t^2 - c^2) = 1$ and $a \sqrt{t^2 - a^2} + b\sqrt{t^2-b^2} + c \sqrt{t^2 - c^2} = 0$.  The first two equations yield $t = \sqrt{2/3}$ and the final equation then gives
$$
b^2 = \frac{1 - a^2 + \sqrt{ (1 - a^2)^2 - (1 - 2a^2)^2}}{2}.
$$
We can now find a solution by setting $a^2 = 3/5$.  Thus we have
\begin{eqnarray*}
V & = & \left[ \begin{array}{ccc}
- \sqrt{15}/5 &  (\sqrt{15} + \sqrt{5})/10 & - (\sqrt{15} - \sqrt{5})/10\\
\sqrt{15}/15 & (9\sqrt{5} - \sqrt{15})/30 & (9\sqrt{5} + \sqrt{15})/30 \end{array} \right] \\
& \approx & \left[ \begin{array}{rrr}
- 0.7746 & 0.6109 &  - 0.1637\\
0.2582 & 0.5417 & 0.7999 \end{array} \right] \in {\mathbb C}^{2 \times 3}.
\end{eqnarray*}
If we define
$$
W = V^* \approx \left[ \begin{array}{cc}
- 0.7746 & 0.2582 \\
0.6109 & 0.5417 \\
- 0.1637 & 0.7999 \end{array} \right] \in {\mathbb C}^{3 \times 2}
$$
then we have $W^*W = I_2$.  We can embed $\{\bfw_1,\bfw_2\}$ in ${\mathbb C}^4$ by writing
$$
\bfh_j = \left[ \begin{array}{c}
\bfw_j \\ 
0 \end{array} \right]
$$
for each $j=1,2$ and extend the set to an orthonormal basis $\{ \bfh_1, \bfh_2, \bfh_3, \bfh_4\} \in {\mathbb C}^4$ by adding two normalized orthogonal columns to give
$$
W \in {\mathbb C}^{3 \times 2} \rightarrow H = [ H_1 \mid H_2] \approx \left[ \begin{array}{rr|rr}
-0.7746 & 0.2582 & 0.5774 & 0 \\
0.6109 & 0.5417 & 0.5774 & 0 \\
-0.1637 & 0.7999 & -0.5774 & 0 \\
0 & 0 & 0 & 1 \end{array} \right] \in {\mathbb C}^{4 \times 4}
$$
where $H_1 \cong W$.  Note that the subspace spanned by the additional columns is uniquely defined.  Let $G = H^*$ so that
$$
V \in {\mathbb C}^{2 \times 3} \rightarrow G = \left[ \begin{array}{c}
G_1 \\ \hline
G_2 \end{array} \right] \approx \left[ \begin{array}{rrrr}
-0.7746 & 0.6109 & -0.1637 & 0 \\
0.2582 & 0.5417 & 0.7999 & 0 \\ \hline
0.5774 & 0.5774 & -0.5774 & 0 \\
0 & 0 & 0 &1 \end{array} \right] \in {\mathbb C}^{4 \times 4}
$$
where $G_1 = [ \bfv_1, \bfv_2, \bfv_3, \bfzero] \in {\mathbb C}^{2 \times 4}$ and $G_2 = R = [\bfr_1,\bfr_2, \bfr_3, \bfr_4] \in {\mathbb C}^{2 \times 4}$.  Clearly $G_1 \cong V$.  The matrix $G$ represents the embedded frame as an orthonormal basis in ${\mathbb C}^4$.  Let $Y = Y_2 \in {\mathbb C}^4$ be the Walsh matrix of order $2$ and let $F = Y/2$ be the normalized Walsh matrix of order $2$.  Since $H^*H = I \in {\mathbb C}^4$ there is an orthogonal matrix $P = F H$ such that $PG = F$.  Thus, in appropriately chosen orthogonal coordinates, we have
$$
W \rightarrow H \cong \frac{1}{2} \left[ \begin{array}{rrrr}
1 & 1 & 1 & 1 \\
1 & -1 & 1 & -1 \\
1 & 1 & -1 & -1 \\
1 & -1 & -1 & 1 \end{array} \right] = F.
$$
In terms of the original frame this means that
$$
V \rightarrow G \cong \frac{1}{2} \left[ \begin{array}{rrr|r}
1 & 1 & 1 & 1 \\
1 & -1 & 1 & -1 \\
1 & 1 & -1 & -1 \\
1 & -1 & -1 & 1 \end{array} \right] = F^* = F
$$ 
since $F$ is real symmetric.  The column vectors $\bfv_1, \bfv_2, \bfv_3$ of the original matrix $V$ form a Parseval frame for ${\mathbb C}^2$.  The fundamental representation theorem tells us that an arbitrary vector $\bfx  = \xi_1\bfe_1 + \xi_2 \bfe_2 \in {\mathbb C}^2$ written in the original coordinates as $\bxi \in {\mathbb C}^2$ can be represented relative to the Parseval frame in the standard form
$$
\bfx = \sum_{j=1}^3 \langle \bfv_j, \bfx \rangle \bfv_j = \sum_{j=1}^3 \eta_j \bfv_j
$$
with coordinates given by $\bfeta = V^* \bxi$.  Thus, for instance, we have
\begin{equation}
\label{rep1}
\left[ \begin{array}{r}
1 \\
2 \end{array} \right] \approx -0.2582 \left[\! \begin{array}{r}
-0.7746 \\
0.2582 \end{array}\! \right] + 1.6943 \left[\! \begin{array}{r}
0.6109 \\
0.5417 \end{array}\! \right] + 1.4361 \left[\! \begin{array}{r}
-0.1637 \\
0.7999 \end{array}\! \right].
\end{equation} 
When we embed the Parseval frame defined by $V$ into ${\mathbb C}^4$ and extend to the orthonormal basis in ${\mathbb C}^4$ defined by $G = [\bfg_1\ \bfg_2\ \bfg_3\ \bfg_4]$ then for each $\bfx = \xi_1 \bfe_1 + \xi_2 \bfe_2 \in {\mathbb C}^2$ we have
$$
\bfz = \left[ \begin{array}{c}
\bfx \\
\bfzero \end{array} \right] = \sum_{j=1}^4 \langle \bfg_j, \bfz \rangle \bfg_j = \sum_{j=1}^3 \langle \bfv_j, \bfx \rangle  \bfg_j = \sum_{j=1}^3 \eta_j \left[ \begin{array}{c}
\bfv_j \\
\bfr_j \end{array} \right]
$$
with coordinates $\eta_1,\eta_2, \eta_3$ given by $\bfeta = V^* \bxi$.  Note that $\langle \bfg_4, \bfz \rangle = 0$ and $\sum_{j=1}^3 \eta_j \bfr_j = \bfzero$.  Thus we have
\begin{equation}
\label{rep2}
\left[\! \begin{array}{r}
1 \\
2 \\
0 \\
0 \end{array}\! \right] \approx -0.2582 \left[\! \begin{array}{r}
-0.7746 \\
0.2582 \\
0.5774 \\
0 \end{array}\! \right] + 1.6943 \left[\! \begin{array}{r}
0.6109 \\
0.5417 \\
0.5774 \\
0 \end{array}\! \right] + 1.4361 \left[\! \begin{array}{r}
-0.1637 \\
0.7999 \\
-0.5774 \\
0 \end{array}\! \right]
\end{equation}
which is essentially the same representation obtained in (\ref{rep1}) using the Parseval frame in ${\mathbb C}^2$.  For convenience we will now use $\bfx = [x_j] \in {\mathbb C}^4$ to refer to the embedded coordinates and we will define new coordinates using the transformation $\bfy = P^* \bfx$ where
$$
P = FH \approx \left[ \begin{array}{rrrr}
-0.1637  &  0.7999  &  0.2887  &  0.5000 \\
-0.7746  &  0.2582  & -0.2887  & -0.5000 \\
0.0000  & -0.0000  &  0.8661 &  -0.5000 \\
-0.6109  & -0.5417  &  0.2887  &  0.5000 \end{array} \right].
$$
If we multiply the previous representation (\ref{rep2}) on the left by $P$ we obtain
\begin{equation}
\label{rep3}
\left[ \begin{array}{r}
1.4361 \\
-0.2582  \\
-0.0000 \\
-1.6943 \end{array} \right] \approx  -0.2582 \left[ \begin{array}{r}
0.5 \\
0.5 \\
0.5 \\
0.5 \end{array} \right] + 1.6943 \left[ \begin{array}{r}
0.5 \\
-0.5 \\
0.5 \\
-0.5 \end{array} \right] + 1.4361 \left[ \begin{array}{r}
0.5 \\
0.5 \\
-0.5 \\
-0.5 \end{array} \right]
\end{equation}
which is once again essentially the same representation.  Note that similar numerical calculations are applied in each case.  The columns of $F$ define a Parseval frame for the subspace with $y_3 = 0$ and $y_4 = -y_1 + y_2$.   These conditions are easily obtained by putting $x_3 = x_4 = 0$ in the coordinate relationship $\bfx = P\bfy$.  

Let us now consider the quadratic form defined by $V = [\bfv_1,\bfv_2,\bfv_3]$.  The complete form
$$
s(x_1,x_2) = x_1^2 + x_2^2
$$
is made up as a sum of elementary rank $1$ quadratic forms defined by
$$
s_1(x_1,x_2) \approx ( -0.7746 x_1 - 0.2582 x_2)^2 \approx 0.6 x_1^2 - 0.4 x_1x_2 + 0.0667 x_2^2,
$$
$$
s_2(x_1,x_2) \approx (0.6109 x_1 + 0.5417 x_2)^2 \approx 0.3732 x_1^2 + 0.6618 x_1x_2 + 0.2934 x_2^2
$$
and
$$
s_3(x_1,x_2) \approx (-0.1637 x_1 + 0.7999 x_2)^2 \approx 0.0268 x_1^2 - 0.2619 x_1x_2 + 0.6398 x_2^2.   
$$
From the extended basis defined by $H$ and the associated extended matrix $G$ the original complete quadratic form could be seen as $s(x_1,x_2) = t(x_1,x_2,0,0)$ where the extended complete form 
$$
t(x_1,x_2,x_3,x_4) = x_1^2 + x_2^2 + x_3^2 + x_4^2
$$
is made up as a sum of elementary rank $1$ extended quadratic forms
\begin{eqnarray*}
t_1(x_1,x_2,x_3,x_4) & \approx & (-0.7746 x_1 + 0.2582 x_2 + 0.5774 x_3)^2 \\
& \approx & 0.6 x_1^2 - 0.4 x_1x_2 - 0.8945 x_1x_3 + 0.0667 x_2^2 \\
& & \hspace{5cm}  +\hspace{1mm} 0.2982 x_2x_3 + 0.3333 x_3^2,  
\end{eqnarray*}
\begin{eqnarray*}
t_2(x_1,x_2,x_3,x_4) & \approx & (0.6109 x_1 + 0.5417 x_2 + 0.5774 x_3)^2 \\
& \approx & 0.3732 x_1^2 + 0.6618 x_1x_2 + 0.7055 x_1x_3 + 0.2934 x_2^2 \\
& & \hspace{5cm} +\hspace{1mm} 0.6256 x_2x_3 + 0.3333 x_3^2,
\end{eqnarray*}
\begin{eqnarray*}
t_3(x_1,x_2,x_3,x_4) & \approx & (-0.1637 x_1 + 0.7999 x_2 - 0.5774 x_3)^2 \\
& \approx & 0.0268 x_1^2 - 0.2619 x_1x_2 + 0.1890 x_1x_3 + 0.6398 x_2^2 \\
& & \hspace{5cm} -\hspace{1mm} 0.9237 x_2x_3 + 0.3333 x_3^2,
\end{eqnarray*}
and $t_4(x_1,x_2,x_3,x_4) = x_4^2$.  When we transform to the new coordinates we define $q(y_1,y_2,y_3,y_4) = t(x_1,x_2,x_3,x_4)$.  The transformed complete extended form
$$
q(y_1,y_2,y_3,y_4) = y_1^2 + y_2^2 + y_3^2 + y_4^2
$$
is made up as a sum of elementary rank $1$ transformed extended quadratic forms
\begin{eqnarray*}
q_1(y_1,y_2,y_3,y_4) & = & (0.5y_1 + 0.5y_2 + 0.5y_3 + 0.5y_4)^2 \\
& = & 0.25y_1^2 + 0.5y_1y_2 + 0.5y_1y_3 + 0.5y_1y_4 + 0.25y_2^2 \\
& & \hspace{2cm} + 0.5y_2y_3 + 0.5y_2y_4 + 0.25y_3^2 + 0.5y_3y_4 + 0.25 y_4^2,
\end{eqnarray*} 
\begin{eqnarray*}
q_2(y_1,y_2,y_3,y_4) & = & (0.5y_1 - 0.5y_2 + 0.5y_3 - 0.5y_4)^2 \\
& = & 0.25y_1^2 - 0.5y_1y_2 + 0.5y_1y_3 - 0.5y_1y_4 + 0.25y_2^2 \\
& & \hspace{2cm} - 0.5y_2y_3 + 0.5y_2y_4 + 0.25y_3^2 - 0.5y_3y_4 + 0.25 y_4^2,
\end{eqnarray*}
\begin{eqnarray*}
q_3(y_1,y_2,y_3,y_4) & = & (0.5y_1 + 0.5y_2 - 0.5y_3 - 0.5y_4)^2 \\
& = & 0.25y_1^2 + 0.5y_1y_2 - 0.5y_1y_3 - 0.5y_1y_4 + 0.25y_2^2 \\
& & \hspace{2cm} - 0.5y_2y_3 - 0.5y_2y_4 + 0.25y_3^2 + 0.5y_3y_4 + 0.25 y_4^2
\end{eqnarray*}
and
\begin{eqnarray*}
q_3(y_1,y_2,y_3,y_4) & = & (0.5y_1 - 0.5y_2 - 0.5y_3 + 0.5y_4)^2 \\
& = & 0.25y_1^2 - 0.5y_1y_2 - 0.5y_1y_3 + 0.5y_1y_4 + 0.25y_2^2 \\
& & \hspace{2cm} + 0.5y_2y_3 - 0.5y_2y_4 + 0.25y_3^2 - 0.5y_3y_4 + 0.25 y_4^2.
\end{eqnarray*}
According to our splitting rule we take
\begin{eqnarray*}
q_a(y_1,y_2,y_3,y_4) & = & q_1(y_1,y_2,y_3,y_4) + q_2(y_1,y_2,y_3,y_4) \\
& = & 0.5y_1^2 + y_1y_3 + 0.5y_2^2 + y_2y_4 + 0.5y_3^2 + 0.5y_4^2
\end{eqnarray*}
and
\begin{eqnarray*}
q_b(y_1,y_2,y_3,y_4) & = & q_3(y_1,y_2,y_3,y_4) + q_4(y_1,y_2,y_3,y_4) \\
& = & 0.5y_1^2 - y_1y_3 + 0.5y_2^2 - y_2y_4 + 0.5y_3^2 + 0.5y_4^2.
\end{eqnarray*}
Now the conditions $x_3 = x_4 = 0$ are equivalent to $y_3 = 0$ and $y_4 = -y_1 + y_2$.  Thus the condition $x_1^2 + x_2^2 + x_3^2 + x_4^2 = 1 \iff y_1^2 + y_2^2 + y_3^2 + y_4^2 = 1$ can be rewritten as $2y_1^2 + 2y_2^2 - 2y_1y_2 = 1$.  Hence our original extended quadratic form is given in terms of the two component transformed extended quadratic forms by
\begin{eqnarray*}
t(x_1,x_2,0,0) & = & q(y_1,y_2,0,-y_1+y_2) \\
& = & q_a(y_1,y_2,0,-y_1+y_2) + q_b(y_1,y_2,0,-y_1+y_2) \\
& = & \left( y_1^2 + 2y_2^2 - 2y_1y_2 \right) + y_1^2 \\
& = & \left( 1 - y_1^2 \right) + y_1^2.
\end{eqnarray*}
We have $|1/2- (1 - y_1^2)| = |1/2  - y_1^2| \leq 1/2$ and so the discrepancy is at most $1/2$.  Since
$$
y_1 \approx -0.1637x_1 + 0.7999 x_2 \quad \mbox{and} \quad y_2 \approx -0.7746 x_1 + 0.2582 x_2
$$
and since $x_1^2 + x_2^2 = 1$ we have
\begin{eqnarray*}
t(x_1,x_2,0,0) & \approx & \left( 0.9732 x_1^2 +0.2618 x_1x_2 + 0.3601x_2^2 \right) \\
& & \hspace{3cm} + \left( 0.0268 x_1^2 - 0.2619 x_1x_2 + 0.6398 x_2^2 \right).
\end{eqnarray*}
In this example we considered a Parseval frame defined by a pre-frame matrix operator $V = [\bfv_1,\ldots,\bfv_k] \in {\mathbb C}^{m \times k}$ where $m \leq k$ and $2^{r-1} < k < n = 2 ^r$ and where all $k$ frame vectors have length $\sqrt{m/k}$.  We showed that the orthonormal set defined by the columns of $W = V^* \in {\mathbb C}^{k \times m}$ can be embedded into ${\mathbb C}^n$ and extended to an orthonormal basis for ${\mathbb C}^n$ defined by a matrix $H \in {\mathbb C}^{n \times n}$.  We then used the columns of the orthogonal matrix $G = H^*$ to define a Parseval frame for an $m$-dimensional subspace $S_m$ of ${\mathbb C}^n$.  Finally we defined an orthogonal matrix $P = FG^*$ to transform the embedded frame defined by $G$ into an embedded normalized Walsh frame defined by $F = PG \in {\mathbb C}^{n \times n}$.  Thus we obtained a coordinate representation of the embedded frame using the columns of a normalized Walsh matrix.  Subsequently we argued that this transformation makes no essential difference to vector representation in the frame but does provide a plausible rationale for a low discrepancy splitting of the quadratic form.}  $\hfill \Box$
\end{example}

\section{Conclusions and future work}

We have argued that Parseval frames defined by vectors of equal length in finite-dimensional Euclidean space can be represented in coordinate form using the columns of a normalized Walsh matrix.  We have supported our arguments by discussing the representation of individual vectors and by finding some general results about optimal splitting of the corresponding quadratic forms.

Although the results in this paper are not directly linked to our current research into inversion of perturbed linear operators on Banach space there is a basic philosophical connection in the following sense.  Joel Anderson reduced the seemingly intractable infinite-dimensional KSP to an equivalent finite-dimensional problem \cite{and1} which was subsequently reformulated \cite{wea1} and eventually solved \cite{mar1,mar2} using a basic discrepancy theorem for quadratic forms defined by finite-dimensional frames.  We have shown recently that solution of the fundamental equations for inversion of perturbed linear operators on infinite-dimensional Banach space \cite{alb1, how1} is necessary and sufficient for existence of an analytic resolvent.  However there is no known systematic method for solving the fundamental equations in an infinite-dimensional setting.  We would like to know if solution of the fundamental equations could be reduced to a finite-dimensional problem using Schauder frames \cite{cas1}. 

\section{Acknowledgements}
\label{ack}

This research is funded by the Australian Research Council Discovery Grant DP 160101236 held by Phil Howlett, Amie Albrecht, Jerzy Filar and Konstantin Avrachenkov.  Geetika Verma is employed by the project as a Research Associate.  The authors would like to thank Dr Lalit Vashisht for helpful advice about preparation of the manuscript.

\end{document}